\documentclass[12pt]{article}

\usepackage{amsmath,amssymb,amsfonts,theorem,makeidx,latexsym,epsfig,,subfigure}

\newtheorem{defn}{Definition}[section]

\newtheorem{lemma}[defn]{Lemma}

{\theorembodyfont{\rmfamily}

\newtheorem{ex}[defn]{Example}}

\newtheorem{thm}[defn]{Theorem}

\newtheorem{prop}[defn]{Proposition}

\newtheorem{cor}[defn]{Corollary}

\numberwithin{equation}{section}

\newcommand{\ltr}{ L^2(\mathbb R) }

\newcommand{\ltz}{{\ell}^2(\mathbb Z)}

\newcommand{\md}{\{ E_{m/M}T_{nN}g \}_{n\in \mz,m=0, \dots,M-1}}

\newcommand{\mdx}{\{ E_{m/M}T_{nN}g \}}

\newcommand{\mn}{\mathbb N}

\newcommand{\mr}{\mathbb R}

\newcommand{\mz}{\mathbb Z}

\def\bp{{\noindent\bf Proof. \ }}

\def\ep{\hfill$\square$\par\bigskip}

\def\bqs{\begin{equation}}

\def\eqs{\tag*{$\square$}\end{equation}\par\bigskip}

\def\nl{\left|\left|}

\def\nr{\right|\right|}

\def\Span{\text{span}}

\def\supp{\text{supp}}

\def\bop{\begin{op}\rm}

\def\eop{\end{op}}

\def\bee{\begin{eqnarray}}

\def\ene{\end{eqnarray}}

\def\bes{\begin{eqnarray*}}

\def\ens{\end{eqnarray*}}

\def\bei{\begin{itemize}}

\def\eni{\end{itemize}}

\def\bt{\begin{thm}}

\def\et{\end{thm}}

\def\bc{\begin{cor}}

\def\ec{\end{cor}}

\def\bpr{\begin{prop}}

\def\epr{\end{prop}}

\def\bl{\begin{lemma}}

\def\el{\end{lemma}}

\def\bd{\begin{defn}}

\def\ed{\end{defn}}

\def\bex{\begin{ex}}

\def\enx{\end{ex}}

\def\bfi{\begin{fig}}

\def\efi{\end{fig}}

\title{Gabor frames  in $\ltz$ and linear dependence}

\date{\today}

\author{Ole Christensen, Marzieh Hasannasab}

\begin{document}

\maketitle

\begin{abstract} We  prove that  an overcomplete Gabor frame in $\ltz$ generated by a finitely
supported sequence is always  linearly dependent. This is a particular case of
a general result about linear dependence versus independence for  Gabor systems in
$\ltz$ with modulation parameter $1/M$ and translation parameter
$N$ for some  $M,N\in \mn,$ and generated by a finite sequence $g$ in $\ltz$ with $K$ nonzero entries.

\end{abstract}

\begin{minipage}{120mm}


\end{minipage}

\section{Introduction}
Linear dependence versus linear independence is a well-studied topic in Gabor analysis. In
particular Linnell \cite{Linn} proved that any Gabor system in $\ltr$ generated by a nonzero function and a time-frequency lattice $a\mz \times b\mz$ is  linearly independent, hereby confirming
a conjecture by Heil, Ramanathan and Topiwala \cite{HRT}. The analogous problem based on
time-frequency shifts on a general locally compact abelian group was studied by Kutyniok in \cite{Kutyn2}
and
Gabor systems on finite groups were analyzed in the paper \cite{Lawrence} by Lawrence, Pfander, and Walnut. Results by Jitomirskaya \cite{jito} imply that the conjecture would
fail on $\ltz,$ as explained by
Demeter and Gautam in \cite{Demeter}.

The purpose of this short note is to give a more detailed discussion of frame properties and linear independence
versus linear dependence for Gabor systems in $\ltz.$ In particular we prove that an
overcomplete
Gabor frame in $\ltz$ generated by a finite sequence is always linearly dependent.
Furthermore we collect and apply various
methods for analysis of such frames, e.g., the duality principle, sampling of
Gabor frames for $\ltr,$ and perturbation methods.
For $g\in \ltz$ we denote the $j$th coordinate
by $g(j).$
For $M\in \mn,$ define the {\it modulation
operators} $E_{m/M}, m=0, \dots, M-1,$ acting on $\ltz$ by $E_{m/M}g(j):= e^{2\pi i jm/M}g(j);$ also, define
the {\it translation operators} $T_{n}, \, n\in \mz,$ by $T_{n}g(j)=g(j-n).$ The {\it Gabor system}
generated by a fixed $g\in \ltz$ and some $M,N\in \mn$ is $\md;$ specifically, $E_{m/M}T_{nN}g$ is
the sequence in $\ltz$  whose $j$th coordinate is
\bes E_{m/M}T_{nN}g(j)= e^{2\pi ijm/M}g(j-nN).\ens
In the rest of this note we will write $\mdx$ instead of $\md.$
It is well-known \cite{CV2} that $\mdx$ can only be a frame for $\ltz$ if $N/M\le 1.$ We prove that if $N/M<1,$
such frames can be constructed with windows $g$ having any number $K\ge N$ of nonzero entries; in contrast to
the case of Gabor frames in $\ltr$ these frames are always linearly dependent. Similarly, for
$M=N$ we can construct Riesz bases for $\ltz$ with
windows $g$ having any number $K\ge N$ of nonzero entries; however, for exactly
the same parameter choices there also exist linearly dependent Gabor systems. More generally,
we characterize  the parameters $M,N,K$ for which  the Gabor system is automatically
linearly independent, linear dependent, resp. that both cases can occur
depending on the choice of $g\in \ltz.$

\section{Gabor systems in $\ltz$ } \label{70307g}

For a finitely supported sequence $g\in \ltz,$ let $|\supp \, g |$ denote the number of nonzero
entries of $g.$
For illustrations and concrete examples we will often  use the sequences
$\delta_k\in \ltz, k\in \mz,$ given by
\bes \delta_k(j)= \begin{cases} 1 &\mbox{if} \, j=k, \\
0 &\mbox{if} \, j\neq k. \end{cases} \ens

It was observed already by Lopez  \& Han \cite{LoHa} that for any $M,N\in \mn$
with $N\le M$ there exist frames $\mdx$ for $\ltz$ generated by windows with $N$
nonzero elements. We will need the following extension, characterizing
the existence of Gabor frames
$\mdx$ for $\ltz$ with a given support size $K.$

\bt \label{70307b}
Let $M,N,K\in \mn.$ Then the following hold:
\bei
\item[(i)] There exists a Gabor frame $\mdx$ for $\ltz$ generated by a window
$g$ with $|\supp \,g|=K$ if and only if  $N\le M$ and $K\ge N$.
\item[(ii)] There exists a Riesz sequence $\mdx$ in $\ltz$ generated by a window
$g$ with $|\supp \,g|=K$ if and only if  $N\ge M$ and $K\ge M$.
\eni
\et

\bp For the proof of (i), the necessity of the condition $N\le M$ is obvious.
We will now show that if $K < N$ then $\mdx$ can not be complete in $\ltz.$ We do this by identifying some $k\in\mz$ such that $E_{m/M}T_{nN}g(k)=0$ for all $n\in\mz$ and $m\in\{0,\dots,M-1\}$. Consider $I:=\{1,\dots,N\}$; then, for any $j\in\mz$, there exists exactly one value of $n\in\mz$ such that $j+nN\in I$. Since $g(j)\neq 0$ only occur for $K<N$ values of $j$, there exists some $k\in I$ such that $j+nN\neq k$ for all $n\in\mz$ and all $j\in\mz$ such that $g(j)\neq 0$. That is, $k-nN\neq j$ for all $n\in\mz$ and all $j\in\mz$ such that $g(j)\neq 0 $. Thus for all $n\in\mz$, we have that $g(k-nN)=0$. This proves that  $E_{m/M}T_{nN}g(k)=0$ for all $n\in\mz$ and $m\in\{0,\dots,M-1\}$ and thus $\mdx$ can not be complete if $K<N$; in other words, $K\geq N$ is necessary for $\mdx$ to be a frame for $\ltz$.

Now assume that $N\le M$ and consider any $g\in \ltz$ for which
\bee \label{70307d} g(j)\neq 0 \, \mbox{for} \ j\in \{1, \dots, N\} \mbox{ and } g(j)=0 \mbox { for } j\notin \{1, \dots, N\}.\ene
All the vectors in $\{E_{m/M}g\}_{m=0, \dots, M-1}$ have support in $\{1, \dots, N\}.$ Writing the
coordinates for these vectors for $j\in \{1, \dots, N\}$ as rows in an $M \times N$ matrix, we get
\bes {\cal A}\ & = & \left(
\begin{array}{ccccc}
g(1) & g(2)  &  \cdot & \cdot  & g(N) \\
e^{\frac{2\pi i}{M} } g(1) & e^{\frac{2\pi i}{M}2 } g(2) &  \cdot & \cdot &
e^{\frac{2\pi i}{M} N} g(N) \\
e^{\frac{2\pi i}{M}2 } g(1) & e^{\frac{2\pi i}{M}2\cdot 2 } g(2) & \cdot & \cdot &
e^{\frac{2\pi i}{M}2 \cdot N } g(N)\\
\cdot & \cdot & \cdot & \cdot  & \cdot \\
\cdot  & \cdot & \cdot & \cdot & \cdot  \\
e^{\frac{2\pi i}{M}(M-1) } g(1) & e^{\frac{2\pi i}{M}(M-1)2 } g(2)  &  \cdot & \cdot & e^{\frac{2\pi i}{M}(M-1)N } g(N)
\end{array}
\right).
\ens Thus, letting $\omega:= e^{\frac{2\pi i}{M}},$
\bee \label{70503a} {\cal A} = \left[ w^{(k-1)j}\right]_{k=1, \dots, M, j=1, \dots, N} \mbox{Diag}(g(1), \dots, g(N)).\ene

Proposition 1.4.3 in \cite{CB} shows that the rows in the matrix ${\cal A}$ form a frame for
$\Span \{\delta_k\}_{k=1}^N$ if and only if the columns in ${\cal A}$ are linearly independent;
since $g(j)\neq 0$ for $j=1, \dots, N$ the linear independence
of the columns  follows from \eqref{70503a}. Applying the translation operators $T_{nN}$
it now follows that $\{E_{m/M}T_{nN}g\}_{n\in \mz, m=0, \dots,M-1}$ is a frame for $\ltz$, with $K=N$.

Now, consider any $K>N$
and any $\epsilon >0$ and
let $\widetilde{g}:= g+ \epsilon \sum_{k=N+1}^K \delta_k.$ It is easy to see that
$\{E_{m/M}T_{nN}\delta_k\}$ is a Bessel sequence with bound $M;$ it follows that
for any finite
sequence $\{c_{m,n}\} \in \ell^2(\{1, \dots, M-1\} \times \mz),$
\bes \nl \sum c_{m,n} E_{m/M}T_{nN} ( \widetilde{g}-g) \nr & = &
\nl \epsilon \sum_{k=N+1}^K\sum c_{m,n} E_{m/M}T_{nN} \delta_k \nr \\
& \le & \epsilon \sum_{k=N+1}^K \nl \sum c_{m,n} E_{m/M}T_{nN} \delta_k \nr \\
& \le & \epsilon (K-N)\sqrt{M}  \left(\sum | c_{m,n}|^2 \right)^{1/2}.\ens
Let $A$ denote a lower frame bound for $\{E_{m/M}T_{nN}g\}_{n\in\mz, m=0,\ldots,M-1}.$ If we choose $\epsilon >0$ such that
$\epsilon (K-N)\sqrt{M} < A,$ it follows from Theorem 22.1.1 in \cite{CB} that $\{E_{m/M}T_{nN}\widetilde{g}\}_{m=0, \dots, M-1, n\in \mz}$
is a frame for $\ltz.$ By construction, $K= | \supp \, g |.$

The result in (ii) is  a consequence of the duality principle \cite{LeSi2}, stating that
a Bessel sequence $\mdx$ is a frame for $\ltz$ if and only if the Gabor
system $\{E_{m/N}T_{nM}g\}$ is a Riesz sequence; in particular
the finitely
supported windows $g$ generating frames in (i) are precisely the ones that generate
Riesz sequences in (ii).
A direct proof of the existence can be given along the lines of the proof of (i), as
follows.
Assume that $M\le N$ and consider any $g\in \ltz$ for which
$g(j)\neq 0 \, \mbox{for} \ j\in \{1, \dots, M\} \mbox{ and } g(j)=0 \mbox { for } j\notin \{1, \dots, M\}.$ Then $\{E_{m/M}g\}_{m=0, \dots, M-1}$ is a basis for
$\Span \{\delta_k\}_{k=1}^M;$  since $N \ge M$ this implies that $\mdx$ is a Riesz sequence
in $\ltz.$  A similar perturbation argument as in (i) now yields the conclusion.  \ep

Let us mention yet another way
of proving the existence of Gabor frames
$\mdx$ for $N/M<1,$ using sampling of B-spline generated Gabor frames for $\ltr.$
Recall that the B-splines $B_K, K\in \mn,$ are defined recursively by convolutions,
$B_1:= \chi_{[0,1]}, B_{K+1}(x):= (B_K * B_1)(x)= \int_0^1 B_K(x-t)\, dt, \, x\in \mr.$

\bex Assume that $N<M$ and consider the B-spline $B_{N+1}$. Since
    $1/M\leq 1/(N+1)$, the system $ \{e^{2\pi i mx/M}B_{N+1}(x -nN)\}_{n,m\in\mz} $ is a Gabor frame for $L^2(\mr)$ by Corollary 11.7.1 in \cite{CB}. Define the discrete sequence $B_{N+1}^D=\{B_{N+1}(j)\}_{j\in\mz}$. Since $B_{N+1}$ is a continuous function with compact support, the sampling results in
    \cite{Jan6} imply that the discrete Gabor system $\{E_{m/M}T_{nN}B_{N+1}^D\}_{n\in\mz, m=0,\ldots,M-1} $ is a frame for $\ltz$. Note that
     $\supp \, B_{N+1}^D= \{1, 2, \dots, N\},$ i.e., $|\supp~ B_{N+1}^D|=N$.
\ep \enx

The main body of Gabor analysis in $\ltr$ has a completely parallel version
in $\ltz,$ but with regard to linear dependence the two cases are
very different. In fact, certain choices of the parameters $M,N,K\in \mn$
imply that the Gabor system  $\mdx$ is linearly dependent for all
windows $g\in \ltz$ with $|\supp \, g |=K;$ for other choices of the
parameters there exist linearly dependent as well as linearly independent
Gabor systems.  The precise statement is as follows.

\bt \label{70307a} Let $M,N\in \mn.$ Then the following hold:
\bei
\item[(i)] If $M=1,$ the system $\mdx$ is linearly independent for all $g\in \ltz \setminus \{ 0\}.$
\item[(ii)] If $M>|\supp \, g |$ the Gabor system $\mdx$ is linearly dependent.
\item[(iii)] If $N<M,$ the Gabor system $\mdx$ is linearly dependent
for any finitely supported $g\in \ltz$.
\item[(iv)] For all $M,N,K\in\mn$ there exists a linearly
dependent Gabor system $\mdx$ with $K= |\supp \, g |.$
\item[(v)] If $N\ge M,$ then there exists for any $K\ge M$ a linearly
independent Gabor system $\mdx$ with $K= |\supp \, g |.$
\eni
\et

\bp For $M=1$ the system $\mdx$ equals the shift-invariant system $\{T_{nN}g\}_{n\in \mz}$
and is thus linearly independent whenever $g\in \ltz \setminus\{ 0\};$ this proves (i).
For the proof of (ii),  the vectors $\{E_{m/M}g\}_{m=1, \dots, M-1}$ can be considered as $M$
vectors in a space of dimension $|\supp \, g|;$ thus they are linearly dependent if $M>|\supp \, g|,$
and hence $\mdx$ is linearly dependent.

For the proof of (iii), consider any finitely supported $g\in \ltz.$ Without loss of generality, assume that
$g(j)=0$ for $j\notin \{1,2, \dots L\}.$ Now, if $L<M$, then the finite collection of vectors  $\{E_{m/M}g\}_{m=0,\ldots,M-1}$ is clearly linear dependent. Thus, we now consider the case $M \leq L $. Considering a finite number of translates of $g$, i.e., $\{T_{n N}g\}_{n=0,\ldots,\ell}$ for some $\ell\in \mn,$ there are at most
$L+ \ell N$ coordinates where one or more of the vectors are nonzero; thus the
 system $\{T_{n N}g\}_{n=0,\ldots,\ell}$ belongs to an $(L+ \ell N)$-dimensional space. Therefore the collection $\{E_{m/M} T_{n N}g\}_{m=0,\ldots,M-1,n=0,\ldots,\ell}$ consists of $( \ell+1)M$ vectors in an $(L+ \ell N)$-dimensional space.  Clearly they are linearly dependent if we choose $\ell\in\mn$ such that $(\ell+1)M >L+ \ell N$, i.e.,  $\ell>\frac{L-M}{M-N}$.
Thus the Gabor system $\mdx$ is linearly dependent, as claimed.

For the proof of (iv), given $M\in \mn,$ let
$g:= \sum_{k=1}^K \delta_{kM};$ then for any $m^\prime\in \mn,$
\bes E_{m^\prime/M}g(j)=  e^{2\pi i m^\prime j/M} \sum_{k=1}^K \delta_{kM}(j)=
\sum_{k=1}^K \delta_{kM}(j)=g(j), \, \forall j\in \mz,\ens i.e., $E_{m^\prime/M}g=g;$ thus
the Gabor system $\mdx$ is linearly dependent.
The result in (v) is a consequence of Theorem \ref{70307b} (ii). \ep

Let us single out the particular result that indeed motivated us to write this
short note. Recall that a frame that is not a basis is said to be {\it overcomplete;}
for a frame $\mdx$ in $\ltz$ this is the case if and only if $N<M$ \cite{CV2}.

\bc  Any overcomplete Gabor frame $\mdx$ with a finitely supported window $g\in \ltz$
is linearly dependent.\ec

\bp The result follows immediately from Theorem \ref{70307a} (iii). \ep

The picture changes if we allow
 windows with infinite support: linearly independent and overcomplete
Gabor frames with infinitely supported windows exist, as we show now.
Our construction is inspired by a calculation for Hermite
functions in $\ltr$ given in \cite{HRT}.

\begin{prop}
      	Define $g\in\ltz$ by $g(j)=e^{-j^2}$.  Then $\mdx$ is   linearly independent for all
      $M,N \in \mn$ and a frame for $\ltz$ if $N<M$.
      \end{prop}
      \bp It is well-known that a Gabor system $\{e^{2\pi ibx} \varphi (x-na)\}_{m,n\in \mz}$
      in $\ltr$ is a Gabor frame for $\ltr$ whenever $\varphi(x)=e^{-x^2}$ and
      $0<ab<1.$ Applying the sampling results by Janssen (see Proposition 2 in \cite{Jan6})
      it follows that the sequence $g$ generates a Gabor frame $\mdx$ for $\ltz$
      whenever $N/M<1.$ Note that this argument uses that the Gaussian satisfies
      the so-called condition R; we refer to \cite{Jan6} for details.

      Now consider any $M,N\in \mn.$ In order to show that $\mdx$ is linearly
      independent,
      assume that there is a finite
      scalar sequence $\{c_{n,m}\}_{n=-L\ldots,L,m=0,\ldots,M-1}$
       such that
      $\sum_{n=-L}^{L}\sum_{m=0}^{M-1} c_{n,m}E_{m/M}T_{nN}g=0.$
      Thus, for all $j\in\mz$,
      \bes\label{1003a}
      0&=&\sum_{n=-L}^{L}\sum_{m=0}^{M-1} c_{n,m} e^{2\pi ijm/M}e^{-(j-nN)^2}
      =e^{-j^2}\sum_{n=-L}^{L}(\sum_{m=0}^{M-1} c_{n,m} e^{2\pi ijm/M})e^{2nNj-(nN)^2}
      \ens For $n=-L, \dots, L,$ defining
      the functions ${\cal E}_n$  on $\mz$
      by ${\cal E}_n(j)=\sum_{m=0}^{M-1} c_{n,m} e^{2\pi ijm/M}, \\ j\in \mz,$ we thus have
      \bee\label{1118a} \sum_{n=-L}^{L}{\cal E}_n(j)e^{2nNj-(nN)^2}=0, \, \forall j\in \mz. \ene
       Note that ${\cal E}_n$ is a bounded and $M$-periodic function on $\ltz$. We will first prove that ${\cal E}_n= 0$ for all $n=-L,\ldots,L$.  Assume that there is some $n>0$ such that ${\cal E}_n(j)\neq0$ for some $j\in\mz$. Then take the largest such $n$ and a corresponding $j_0\in\{1,\ldots,M-1\}$ such that ${\cal E}_n(j_0)\neq 0$. Then
      \[\sum_{n=-L}^{L}{\cal E}_n(j_0+\ell M)e^{-(nN)^2}e^{2nN(j_0+\ell M)}\rightarrow\infty\quad\mbox{ as }\ell\rightarrow\infty\]
      which is contradicting  \eqref{1003a}. Therefore for all $0<n\leq L$, ${\cal E}_n=0$. A similar argument shows that for all $-L\leq n<0$, we have ${\cal E}_n=0$. Now \eqref{1118a} implies that also ${\cal E}_0=0$, as claimed.

      Considering now any $n=-L\ldots,L$, we thus have
      $\sum_{m=0}^{M-1} c_{n,m} e^{2\pi ijm/M}=0$ for all $j=0,\ldots,M-1$.
      Writing this set of equations in matrix form, the matrix describing
      the system is a Vandermonde matrix and thus invertible;
       it follows that $c_{n,m}=0$ for $m=0,\ldots,M-1$. Since $n\in \{-L, \dots, L\}$ was arbitrary, this proves that the Gabor system is linearly independent. \ep

Let us also give a construction of a  linearly dependent Gabor frame
 for $\ltz$ with an infinitely supported window.

\bex   Assume that $N<M$ and consider the sequence
$g\in \ltz$ given by $g(j)=1$ for $j\in \{1, \dots, N\}$ and
 $g(j)=0$ for $ j\notin \{1, \dots, N\}.$
As we have seen in the proof of Theorem \ref{70307b} (i), the system $\mdx$ is a frame
for $\ltz$. For $\epsilon>0$, let $\widetilde{g}=g+\sum_{\ell=1}^{\infty}\frac{\epsilon}{2^\ell}\delta_{\ell M+1}$. Then $\widetilde{g}$ has infinite support and
a similar calculation as in the proof of Theorem \ref{70307b} (i) shows that
for any finite sequence $\{c_{m,n}\}$,
$\|\sum c_{m,n}E_{m/M}T_{nN}(g-\widetilde{g})\|\le \epsilon \sqrt{M}(\sum|c_{m,n}|^2)^{1/2}.$
Applying again the perturbation results for frames (Theorem 22.1.1 in \cite{CB}), it follows that
for sufficiently small $\epsilon$, the system $\{E_{m/M}T_{nN}\widetilde{g}\}$ is a frame for $\ltz$. Now, since $N<M$ and the support of $g$ has length $N,$  the system $\{E_{m/M}g\}_{m=0, \dots, M-1}$ is
      linearly dependent; thus, we can choose a nonzero scalar sequence $\{c_m\}_{m=0}^{M-1}$
      such that $\sum_{m=0}^{M-1}  c_{m}E_{m/M}g=0,$ i.e.,
$\sum_{m=0}^{M-1}  c_{m} e^{2\pi ijm/M}=0$ for $j=1, \dots, N.$ It follows that for
any $\ell\in \mn,$
\bes \sum_{m=0}^{M-1}  c_{m}E_{m/M}\delta_{\ell M+1}(\ell M+1)=
\sum_{m=0}^{M-1}  c_{m}e^{2\pi i(\ell M+1)m/M}=\sum_{m=0}^{M-1}  c_{m}e^{2\pi im/M}=0,\ens
and thus $\sum_{m=0}^{M-1}  c_{m}E_{m/M}\delta_{\ell M+1}=0.$ The construction of
the sequence $\widetilde{g}$ now shows that $ \sum_{m=0}^{M-1}  c_{m}E_{m/M}\widetilde{g}=0;$ it follows that
the Gabor system $\{E_{m/M}T_{nN}\widetilde{g}\}$ is linearly dependent, as claimed. \ep \enx

%We can now complete the argument by showing  that
%      \bee \label{70916a} \sum_{m=0}^{M-1}  c_{m}E_{m/M}g=0 \Leftrightarrow  \sum_{m=0}^{M-1}  %c_{m}E_{m/M}\widetilde{g}=0;\ene indeed, since $\{E_{m/M}g\}_{m=0, \dots, M-1}$ is
%     linearly dependent, it follows from \eqref{70916a}
%      that the Gabor system $\{E_{m/M}T_{nN}\widetilde{g}\}$ is linearly dependent.

\noindent{\bf Acknowledgment:} The authors would like to thank Guido Janssen, Chris Heil and Shahaf
Nitzan for useful comments and references.

\begin{tabbing}
text-text-text-text-text-text-text-text-text-text \= text \kill \\
Ole Christensen \> Marzieh Hasannasab \\
Technical University of Denmark \> Technical University of Denmark  \\
DTU Compute \> DTU Compute \\
Building 303, 2800 Lyngby \> Building 303, 2800 Lyngby \\
Denmark \> Denmark \\
Email: ochr@dtu.dk \> mhas@dtu.dk
\end{tabbing}

\end{document}